\documentclass[12pt]{article}

\usepackage{amsfonts}
\usepackage{amsmath,amsthm}
\usepackage{amssymb}
\usepackage{cite}

\usepackage[unicode,colorlinks,plainpages=false]{hyperref}

\newtheorem{theorem}{Theorem}%  meant for continuous numbers
\newtheorem{lemma}[theorem]{Lemma}
\newtheorem{corollary}[theorem]{Corollary}
\theoremstyle{thmstyletwo}%

\begin{document}

\centerline{\bf A subdirect decomposition of a semigroup of all fuzzy sets}
\centerline{\bf in a semigroup} 

\bigskip

\centerline{Attila Nagy \footnote{email: nagyat@math.bme.hu\\{\bf Keywords}: fuzzy set, ideal, congruence, semigroup, {\bf MSC Classification}: 20M10, 08A72, 
This work was supported by the National Research, Development and Innovation Office – NKFIH, 115288}}

\medskip

\centerline{Department of Algebra, Institute of Mathematics}
\centerline{Budapest University of Technology and Economics}
\centerline{M\H uegyetem rkp. 3., H-1111 Budapest, Hungary}

\begin{abstract}In this paper we give a subdirect decomposition of semigroups $({\mathfrak F}(S); \circ )$, where $S$ is a semigroup, ${\mathfrak F}(S)$ is the set of all fuzzy sets in $S$, and the operation $\circ$ on ${\mathfrak F}(S)$ is defined by the following way: for
$f, g\in {\mathfrak F}(S)$ and $s\in S$, $(f\circ g)(s)=\vee _{{x, y\in S}\atop{s=xy}}(f(x)\wedge g(y))$ if $s\in S^2$, and $(f\circ g)(s)=0$ otherwise.
\end{abstract}

\section{Introduction and motivation}

A function from a nonempty set $S$ into the closed interval $[0, 1]$ of real numbers is called a fuzzy set in $S$. The concept of the fuzzy set was introduced by L.A. Zadeh in \cite{Zadeh}, and has been applied by many authors to the elementary theory of groupoids, semigroups and groups; see, for example, \cite{Dib, Kehayopulu, Kehayopulu2, Kim, Kuroki1, Kuroki2, Kuroki3, Kuroki4}, \cite{Wang, McLean, Mordeson, Nagy5}, and \cite{Rosenfeld}.
By \cite[Proposition 1.3]{Wang} or \cite[Section 2.3]{Mordeson}, if $S$ is a semigroup (written multiplicatively), then the set ${\mathfrak F}(S)$ of all fuzzy sets in $S$ is also a semigroup under the operation $\circ$ defined by the following way: for arbitrary $f, g\in {\mathfrak F}(S)$ and $s\in S$,
\[(f\circ g)(s)=\begin{cases}
\vee _{{x, y\in S}\atop{s=xy}}(f(x)\wedge g(y)), & \text{if $s\in S^2$}\\
0, & \text{otherwise.}\end{cases}\]
In \cite{Nagy5}, the structure of the semigroup $({\mathfrak F}(S);\circ )$ is examined in that case when $S$ is a right zero semigroup.
It is proved in \cite[Theorem 2.6]{Nagy5} that if $S$ is a right zero semigroup, then the semigroup $({\mathfrak F}(S);\circ )$ is a right regular band.
The purpose of our present paper is to examine the structure of the semigroup
$({\mathfrak F}(S);\circ )$, where $S$ is an arbitrary semigroup.
For every nonempty subset $A$ of a semigroup $S$ having the property that $S\setminus A$ is either empty or an ideal of $S$, we define an associative operation $\star$ on the set ${\mathfrak F}(A)$ of all fuzzy sets in $A$, and show that, if $S$ is a union of its subsets $A_i$, $i\in I$ having the above property, then the semigroup $({\mathfrak F}(S);\circ )$ is a subdirect product of semigroups $({\mathfrak F}(A_i);\star )$, $i\in I$. Applying this result, we examine the effect of the divisibility between elements in a semigroup $S$ on the structure of the semigroup $({\mathfrak F}(S); \circ )$. We show that, for an arbitrary semigroup $S$, the semigroup $({\mathfrak F}(S);\circ )$ is a subdirect product of semigroups $({\mathfrak F}(D_a);\star )$, $a\in S$, where $D_a$ denotes the set of all divisors of an element $a$ of $S$.

\medskip

For notions and notations not defined but used in this paper, we refer to books \cite{Clifford1}, \cite{Nagy}, and \cite{Petrich}.

\section{A subdirect decomposition of $({\mathfrak F}(S); \circ )$}

A nonempty subset $J$ of a semigroup $S$ is an \emph{ideal} of $S$ if $SJ\subseteq J$ and $JS\subseteq J$.
The first theorem of this section shows
that we can define an associative operation on the set ${\mathfrak F}(A)$ if A is a nonempty subset of a semigroup $S$ with the property that $S\setminus A$ is either empty or an ideal of $S$.

\begin{theorem}\label{thm1} Let $S$ be a semigroup and $A$ be a nonempty subset of $S$ such that $S\setminus A$ is either empty or an ideal of $S$. For every $f, g\in {\mathfrak F}(A)$, let $f\star g$ be the fuzzy set in $A$ defined by the following way: for arbitrary $a\in A$, let
\[(f\star g)(a)=\begin{cases}
\vee _{{x, y\in S}\atop{a=xy}}(f(x)\wedge g(y)), & \text{if $a\in S^2$}\\
0, & \text{otherwise.}\end{cases}\]
Then $\star$ is an associative operation on ${\mathfrak F}(A)$.
\end{theorem}

\begin{proof} First of all, we note that $\star$ is a well defined operation on ${\mathfrak F}(A)$, because it follows from the assumption for $A$ that if $xy\in A$ ($x, y\in S$), then $x, y\in A$.
To show that the operation $\star$ associative, let $f, g, h$ be arbitrary fuzzy sets in $A$.
Let $a\in A$ be an arbitrary element. Then
\[((f\star g)\star h)(a)=\begin{cases}
\vee _{{x, y\in S}\atop{a=xy}}((f\star g)(x)\wedge h(y)), & \text{if $a\in S^2$}\\
0, & \text{otherwise.}\end{cases}\]
As \[(f\star g)(x)=\begin{cases}
\vee _{{u, v\in S}\atop{x=uv}}(f(u)\wedge g(v)), & \text{if $x\in S^2$}\\
0, & \text{otherwise,}\end{cases}\]
we have
\[((f\star g)\star h)(a)=\begin{cases}
\vee _{{u, v, y\in S}\atop{a=uvy}}(f(u)\wedge g(v)\wedge h(y)), & \text{if $a\in S^3$}\\
0, & \text{otherwise,}\end{cases}\]
using also the fact that the unit interval $[0, 1]$, ordered in the natural way, is a completely distributive complete lattice (\cite{Raney}, \cite{Szasz}), and so
\[(\vee _{i\in I} a_i)\wedge b=\vee _{i\in I}(a_i\wedge b)\] is satisfied for every elements $a_i$ ($i\in I$) and $b$ of $[0, 1]$.
We can prove in a similar way that
\[(f\star (g\star h))(s)=\begin{cases}
\vee _{{u, v, y\in S}\atop{a=uvy}}(f(u)\wedge g(v)\wedge h(y)), & \text{if $a\in S^3$}\\
0, & \text{otherwise.}\end{cases}\]
Thus
\[(f\star g)\star h=f\star (g\star h),\] that is, the operation $\star$ is associative on ${\mathfrak F}(S)$.
\end{proof}

\medskip

In our investigation the following binary relation plays an important role. For a nonempty subset $A$ of a semigroup $S$, let $\Delta _A$ denote the binary relation on the set ${\mathfrak F}(S)$:
\[\Delta _A=\{ (f, g)\in {\mathfrak F}(S)\times {\mathfrak F}(S): f(a)=g(a) \ \hbox{for all}\ a\in A\}.\]

\begin{theorem}\label{thm2} For every non-empty subset $A$ of a semigroup $S$ having the property that $S\setminus A$ is either empty or an ideal of $S$, the binary relation $\Delta _A$ is a congruence on the semigroup $({\mathfrak F}(S); \circ)$.
\end{theorem}

\begin{proof} It is clear that $\Delta _A$ is an equivalence relation on the set ${\mathfrak F}(S)$.
To show that $\Delta _A$ is a congruence on the semigroup $({\mathfrak F}(S); \circ )$, assume
\[(f_1, g_1)\in \Delta _A\quad \hbox{and}\quad (f_2, g_2)\in \Delta _A\] for elements $f_1, f_2, g_1, g_2\in {\mathfrak F}(S)$. Then
\[f_1(a)=g_1(a)\quad \hbox{and}\quad f_2(a)=g_2(a)\] for all $a\in A$.
We show that $(f_1\circ f_2, g_1\circ g_2)\in \Delta _A$.
Let $a\in A$ be an arbitrary element.
If $a\notin S^2$ then
\[(f_1\circ f_2)(a)=0=(g_1\circ g_2)(a).\]
Consider the case when $a\in S^2$. If $a=xy$ for some $x, y \in S$, then $x, y \in A$, because $S\setminus A$ is either empty or an ideal of $S$. Thus
\[f_1(x)=g_1(x)\quad \hbox{and}\quad f_2(y)=g_2(y)\]
from which it follows that
\[(f_1\circ f_2)(a)=\]
\[=\vee _{{x, y \in S}\atop{a=xy}}(f_1(x)\wedge f_2(y))=\vee _{{x, y \in S}\atop{a=xy}}(g_1(x)\wedge g_2(y))=\]
\[=(g_1\circ g_2)(a).\]
Hence \[(f_1\circ f_2, g_1\circ g_2)\in \Delta _A,\]
and consequently $\Delta _A$ is a congruence on the semigroups $({\mathfrak F}(S); \circ )$.
\end{proof}

\medskip
The following theorem is about the connection between the factor semigroup $({\mathfrak F}(S);\circ )/\Delta _A$ and the semigroup
$({\mathfrak F}(A);\star )$.

\begin{theorem}\label{thm3} If $A$ is a nonempty subset of a semigroup $S$ such that $S\setminus A$ is either empty or an ideal of $S$, then the factor semigroup $({\mathfrak F}(S);\circ )/\Delta _A$ is isomorphic to the semigroup $({\mathfrak F}(A);\star )$.
\end{theorem}

\begin{proof} Let $[f]_{\Delta _ A}$ denote the $\Delta _A$-class of the semigroup $({\mathfrak F}(S); \circ )$ containing the fuzzy set $f$ in $S$.
Let $\Psi$ denote the mapping of the factor semigroup $({\mathfrak F}(S); \circ )/\Delta _A$ into the semigroup $({\mathfrak F}(A); \star )$ defined by the following way. For an arbitrary element $[f]_{\Delta _A}\in ({\mathfrak F}(S); \circ )/\Delta _A$, let
\[\Psi ([f]_{\Delta _A})=f^*,\] where $f^*$ is the restriction of $f$ to $A$. Since the equation $[f_1]_{\Delta _A}=[f_2]_{\Delta _A}$ is equivalent to the equation $f^*_1=f^*_2$, then $\Psi$ is a well defined injective mapping of the factor semigroup $({\mathfrak F}(S); \circ )/\Delta _A$ into the semigroup $({\mathfrak F}(A); \star )$. It is clear that $\Psi$ is also surjective. It remains to show that $\Psi$ is a homomorphism. Let the operation on the factor semigroup $({\mathfrak F}(S); \circ )/\Delta _A$ also be denoted by $\circ$. Let \[[f_1]_{\Delta _A}, [f_2]_{\Delta _A}\in ({\mathfrak F}(S); \circ )/\Delta _A\] be arbitrary elements. Then
\[\Psi([f_1]_{\Delta _A}\circ [f_2]_{\Delta _A})=\Psi([f_1\circ f_2]_{\Delta _A})=(f_1\circ f_2)^*\]
and
\[\Psi([f_1]_{\Delta _A})\star \Psi([f_2]_{\Delta _A})=f^*_1\star f^*_2.\]
We show that $(f_1\circ f_2)^*=f^*_1\star f^*_2$.
Let $a\in A$ be an arbitrary element.

\noindent
If $a\notin S^2$, then
\[(f_1\circ f_2)^*(a)=(f_1\circ f_2)(a)=0=(f^*_1\star f^*_2)(a).\]

\noindent
If $a\in S^2$, then $a=xy$ ($x, y\in S$) implies $x, y\in A$, and hence
\[(f_1\circ f_2)^*(a)=\]
\[=(f_1\circ f_2)(a)=\vee _{{x, y\in S}\atop{a=xy}}(f_1(x)\wedge f_2(y))=\vee _{{x, y\in S}\atop{a=xy}}(f^*_1(x)\wedge f^*_2(y))=\]
\[=(f^*_1\star f^*_2)(a).\]
Thus \[(f_1\circ f_2)^*=f^*_1\star f^*_2.\] Consequently
\[\Psi([f_1]_{\Delta _A}\circ [f_2]_{\Delta _A})=\Psi([f_1]_{\Delta _A})\star \Psi([f_2]_{\Delta _A}).\]
Then $\Psi$ is a homomorphism.
\end{proof}

\medskip
The concept of the subdirect product of semigroups is a very important tool of the theory of semigroups. Many papers are published on the subdirect product of semigroups; see, for example, \cite{Chrislock}, \cite{Ciric}, \cite{Kehayopulu0}, \cite{Mitsch}, \cite{Nagy1, Nagy2, Nagy3}, \cite{Nagy4}, \cite{Petrich0}, \cite{Schein}, and \cite{Yamada}. We say that the semigroup $S$ is a \emph{subdirect product} of semigroups $S_i$, $i\in I$, if $S$ is isomorphic to a subsemigroup $T$ of the direct product of semigroups $S_i$, $i\in I$ such that the restriction of the projection homomorphisms to $T$ are surjective. Subdirect decompositions of a semigroup $S$ are closely connected with congruences on $S$.
If $\alpha _i$, $i\in I$ are congruences on a semigroup $S$ such that $\bigcap _{i\in I}\alpha _i$ is the identity relation on $S$, then $S$ is a subdirect product of the factor semigroups $S/\alpha _i$, $i\in I$, and conversely, if a semigroup $S$ is a subdirect product of semigroups $S_i$, $i\in I$, then there are congruences $\alpha _i$, $i\in I$ on $S$ such that $\bigcap_{i\in I}\alpha _i$ is the identity relation on $S$ and $S/\alpha _i\cong S_i$ for every $i\in I$ (see, for example, \cite[I.3.6]{Petrich}).

\medskip
The next theorem gives a subdirect decomposition of the semigroup $({\mathfrak F}(S);\circ )$.

\begin{theorem}\label{thm4} Let $A_i$, $i\in I$ be a nonempty family of nonempty subsets $A_i$ of a semigroup $S$ having the property that
$S\setminus A_i$ is either empty or an ideal of $S$. If $S=\bigcup_{i\in I}A_i$, then the semigroup $({\mathfrak F}(S);\circ )$ is a subdirect product of semigroups $({\mathfrak F}(A_i);\star )$, $i\in I$.
\end{theorem}

\begin{proof} Let $S$ be an arbitrary semigroup. By Theorem~\ref{thm2}, $\Delta _{A_i}$ is a congruence on $({\mathfrak F}(S);\circ)$ for every $i\in I$. Assume $(f, g)\in \bigcap _{i\in I}\Delta _{A_i}$ for
$f, g\in {\mathfrak F}(S)$. Then $(f, g)\in \Delta _{A_i}$ for every $i\in I$. Since $\bigcup _{i\in I}A_i=S$, then we have $f(s)=g(s)$ for every $s\in S$. Thus $f=g$. Consequently
$\bigcap _{i\in I}\Delta _{A_i}=\iota _{{\mathfrak F}(S)}$,
where $\iota _{{\mathfrak F}(S)}$ denotes the identity relation on ${\mathfrak F}(S)$.
Then the semigroup $({\mathfrak F}(S); \circ )$ is a subdirect product of the factor semigroups $({\mathfrak F}(S); \circ )/\Delta _{A_i};\ i\in I$.
By Theorem~\ref{thm1}, $({\mathfrak F}(A_i); \star )$ is a semigroup for every $i\in I$. By Theorem~\ref{thm3}, semigroups $({\mathfrak F}(S); \circ )/\Delta _{A_i}$ and $({\mathfrak F}(A_i); \star )$ are isomorphic for every $i\in I$, and hence the semigroup $({\mathfrak F}(S); \circ )$ is a subdirect product of semigroups $({\mathfrak F}(A_i); \star );\ i\in I$.
\end{proof}

\medskip
In the remainder of this section we apply the previous results to show that the divisors of elements of a semigroup $S$ determine a subdirect decomposition of the semigroup $({\mathfrak F}(S);\circ )$.

\medskip
For an arbitrary element $a$ of a semigroup $S$, let $D_a$ denote the set of all divisors of $a$, that is, $D_a=\{ s\in S: a\in S^1sS^1\}$.

\begin{lemma}\label{lem1} For every element $a$ of a semigroup $S$, $D_a$ is a nonempty subset of $S$ having the property that $S\setminus D_a$ is either empty or an ideal of $S$.
\end{lemma}

\begin{proof} Let $a$ be an arbitrary element of a semigroup $S$. $D_a$ is not empty, because $a\in D_a$. Assume $xy\in D_a$ for elements $x, y\in S$. Then $a\in S^1xyS^1\subseteq S^1xS^1\bigcap S^1yS^1$ from which it follows that $x, y\in D_a$. Hence $S\setminus D_a$ is either empty or an ideal of $S$.
\end{proof}

\begin{theorem}\label{thm11} For an arbitrary semigroup $S$, the semigroup $({\mathfrak F}(S); \circ )$ is a subdirect product of semigroups $({\mathfrak F}(D_a); \star )$, $a\in S$.
\end{theorem}

\begin{proof} Since $a\in D_a$ for every element $a$ of a semigroup $S$, we have $S=\bigcup_{a\in S}D_a$. Thus Lemma~\ref{lem1} and Theorem~\ref{thm4} together imply that the semigroup $({\mathfrak F}(S); \circ )$ is a subdirect product of semigroups $({\mathfrak F}(D_a); \star )$, $a\in S$.
\end{proof}

\bigskip

By Lemma~\ref{lem1} and the proof of Theorem~\ref{thm4},
if $\Delta _{D_{a}}=\iota _{{\mathfrak F}(S)}$ is satisfied for an element $a$ of a semigroup $S$, then ${\mathfrak F}(S)$ is isomorphic to the subdirect factor $({\mathfrak F}(D_a);\star )$. In the next theorem we give a necessary and sufficient condition for $\Delta _{D_a}\neq \iota _{{\mathfrak F}(S)}$ to be satisfied in a semigroup $S$ for all $a\in S$.

\medskip
The least ideal of a semigroup $S$ (if it exists) is called the \emph{kernel} of $S$.

\begin{theorem}\label{thm5} $\Delta _{D_a}=\iota _{{\mathfrak F}(S)}$ is satisfied for an element $a$ of a semigroup $S$ if and only if $a$ is in the kernel of $S$. Thus $\Delta _{D_a}\neq \iota _{{\mathfrak F}(S)}$ is satisfied for all elements $a$ of a semigroup $S$ if and only if $S$ has no a kernel.
\end{theorem}

\begin{proof} Let $a$ be an arbitrary element of a semigroup $S$. $\Delta _{D_a}=\iota _{{\mathfrak F}(S)}$ is satisfied if and only if $D_a=S$, that is, $a\in S^1sS^1$ for every $s\in S$. As $S^1sS^1$ is an ideal of $S$, this last condition is equivalent to the condition that $a$ belongs to every ideal of $S$, that is, $a$ is in the kernel of $S$. The second assertion of the proposition is an obvious consequence of the firs one.
\end{proof}

\section{On the restriction of $\Delta _A$ to $S$}

Let $S$ be an arbitrary semigroup and $A$ be a non empty subset of $S$. Then $C_A : S\mapsto [0, 1]$ defined by
\[C_A(x)=\begin{cases}1, & \text{if $x\in A$}\\
0, & \text{otherwise}\end{cases}\]
is a fuzzy set in $S$; this is the {\it characteristic function} of the subset $A$.
If $s$ is an arbitrary element of a semigroup $S$, then let $s$ denote the one-element subset of $S$ containing the element $s$.
For an arbitrary semigroup $S$,
\[\Phi : s\mapsto C_s\quad (s\in S)\]
is an injective mapping of $S$ into the semigroup $({\mathfrak F}(S); \circ )$.
It is a matter of checking to see that $C_s\circ C_t=C_{st}$ for every $s, t\in S$, and hence $\Phi$ is a homomorphism (see also \cite{Kim}). Hence $\Phi$ is an embedding of $S$ into the semigroup $({\mathfrak F}(S); \circ )$.
In the next part of the paper, we will identify $S$ to the subsemigroup $\Phi (S)$ of the semigroup $({\mathfrak F}(S); \circ )$, where necessary.

\bigskip
\noindent
Let $J$ be a subset of a semigroup $S$ such that $J=\emptyset$ or $J$ is an ideal of $S$. Then
\[\varrho _J=\{ (x, y)\in S\times S:\ x=y \ \hbox{or}\ x,y\in J\}\]
is a congruence on $S$ which is called the \emph{Rees congruence} defined by $J$ (see \cite{Clifford1}). We note that if $J=\emptyset$ then $\varrho _J$ is the identity relation on $S$.

\bigskip
\noindent
For a nonempty subset $A$ of a semigroup $S$ having the property that $S\setminus A$ is either empty or an ideal of $S$, let $\Delta_A\arrowvert_S$ denote the restriction of the congruence $\Delta_A$ on $({\mathfrak F}(S);\circ )$ to $S$ (as a subsemigroup of $({\mathfrak F}(S);\circ )$).

\begin{theorem}\label{thm6} For every nonempty subset $A$ of a semigroup $S$ having the property that $S\setminus A$ is either empty or an ideal of $S$, we have $\Delta _A\arrowvert_S=\varrho _{S\setminus A}$.
\end{theorem}

\begin{proof} Let $S$ be a semigroup and $A$ be a nonempty subset of $S$ having the property that $S\setminus A$ is either empty or an ideal of $S$.
\noindent
First we show that $\Delta_A\arrowvert_S\subseteq \varrho _{S\setminus A}$.
Assume $(s, t)\in \Delta _A$, that is, $(C_s, C_t)\in \Delta _A$ for elements $s, t\in S$ with $s\neq t$. Suppose $s\in A$. Then \[1=C_s(s)=C_t(s)=0\] which is impossible. Similarly, $t\in A$ implies \[0=C_s(t)=C_t(t)=1\] which is also impossible. Hence $s, t\in S\setminus A$, and consequently $(s,t)\in \varrho _{S\setminus A}$.
Thus \[\Delta _A\arrowvert_S\subseteq \varrho _{S\setminus A}.\]

\noindent
To prove the converse inclusion, assume $(s, t)\in \varrho _{S\setminus A}$ for elements $s, t\in S$ with $s\neq t$. Then $s, t\in S\setminus A$.
Let $x\in A$ be an arbitrary element. Then $x\neq s$ and $x\neq t$, and so \[C_s(x)=0=C_t(x).\] Hence $(C_s, C_t)\in \Delta _A$. Thus
\[\varrho _{S\setminus A}\subseteq \Delta _A\arrowvert_S.\]
Consequently $\Delta _A\arrowvert_S=\varrho _{S\setminus A}$.
\end{proof}

\begin{corollary}\label{cor2} For every element $a$ of a semigroup $S$, $\Delta _{D_a}\arrowvert_S=\varrho _{S\setminus D_a}$.
\end{corollary}

\begin{proof} It is obvious by Lemma~\ref{lem1} and Theorem~\ref{thm6}.
\end{proof}

\bigskip

\section{Conclusion} Divisibility between elements in semigroups and subdirect product of semigroups are two important concepts in the study of the structure of semigroups. In this paper we pointed out that the divisibility between elements in a semigroup $S$ determines a subdirect decomposition of the semigroup $({\mathfrak F}(S); \circ )$ of all fuzzy sets in $S$. We proved that, for an arbitrary semigroup $S$, the semigroup $({\mathfrak F}(S); \circ )$ is a subdirect product of semigroups $({\mathfrak F}(D_a); \star )$, $a\in S$, where $D_a$ denotes the set of all divisor of an element $a$ of $S$.


\begin{thebibliography}{130}
\bibitem{Clifford1} Clifford, A.H., Preston, G.B.: The Algebraic Theory of Semigroups I. American Mathematical Society, Providence R.I. (1961)
\bibitem{Chrislock} Chrislock, J.L., Tamura, T.: Notes on subdirect products of semigroups and rectangular bands. Proceedings of the American Mathematical Society 20, 511--514 (1969)
\bibitem{Ciric} \'Ciri\'c, M., Bogdanovi\'c, S.: Subdirect products of a band and a semigroup. Portugaliae Mathematica 53-1, 117--128 (1996)
\bibitem{Dib} Dib, K.A., Galham, N.: Fuzzy ideals and fuzzy bi-ideals in fuzzy semi-groups. Fuzzy Sets and Systems 92, 203-–215 (1997)
\bibitem{Kehayopulu0} Kehayopulu, N., Tsingelis, M.: On subdirectly irreducible ordered semigroups. Semigroup Forum 50, 161-177 (1995)
\bibitem{Kehayopulu} Kehayopulu, N., Tsingelis, M.: Fuzzy sets in ordered groupoids. Semigroup Forum 65, 128--132 (2002)
\bibitem{Kehayopulu2} Kehayopulu, N., Tsingelis, M.: Fuzzy Semiprime and Fuzzy Prime Subsets of Ordered Groupoids. Filomat 31:13 (2017), 4217–-4223
\bibitem{Kim} Kim, K.H.: On fuzzy points in semigroups. International Journal of Mathematics and Mathematical Sciences 26-11, 707--712 (2001)
\bibitem{Kuroki1} Kuroki, N.: Fuzzy bi-ideals in semigroups. Commentarii Mathematici Universitatis Sancti Pauli 28, 17--21 (1979)
\bibitem{Kuroki2} Kuroki, N.: On Fuzzy Ideals and Fuzzy Bi-ideals in Semigroups. Fuzzy Sets and Systems 5, 203--215 (1981)
\bibitem{Kuroki3} Kuroki, N.: On fuzzy semiprime ideals in semigroups. Fuzzy Sets and Systems 8-1, 71-–79 (1982)
\bibitem{Kuroki4} Kuroki, N.: On Fuzzy Semigroups. Information Sciences 53, 203--236 (1991)
\bibitem{Mitsch} Mitsch, H.: Subdirect products of E-inversive semigroups. Journal of the Australian Mathematical Society (Series A) 48, 66--78 (1990)
\bibitem{Wang} Liu, W.: Fuzzy Invariant Subgroups and Fuzzy Ideals. Fuzzy Sets and Systems 8, 133--139 (1982)
\bibitem{McLean} McLean, R.G., Kummer, H.: Fuzzy ideals in semigroups. Fuzzy Sets and Systems 48-1, 137-–140 (1992)
\bibitem{Mordeson} Mordeson, J.N., Malik, D.S., Kuroki, N.: Fuzzy Semigroups. Springer, Berlin (2003)
\bibitem{Nagy1} Nagy, A.: Subdirectly irreducible WE-2 semigroups with globally idempotent core. Lecture Notes in Mathematics 1320, 244--250 (1986)
\bibitem{Nagy2} Nagy, A.: Subdirectly irreducible completely symmetrical semigroups. Semigroup Forum 45, 267--271 (1992)
\bibitem{Nagy3} Nagy, A.: Subdirectly irreducible right commutative semigroups. Semigroup Forum 46, 187--198 (1993)
\bibitem{Nagy} Nagy, A.: Special Classes of Semigroups. Kluwer Academic Publishers, Dordrecht/Boston/London (2001)
\bibitem{Nagy4} Nagy, A.: Subdirectly irreducible RGC-commutative right H-semi\-groups. Semigroup Forum 78-1, 68--76 (2009)
\bibitem{Nagy5} Nagy, A.: Remarks on graphons. International Journal of Algebra 15-2, 61--68 (2021)
\bibitem{Petrich0} Petrich, M.: Regular semigroups which are subdirect products of a band and a semilattice of groups. Glasgow Mathematical Journal 14, 27--49 (1973)
\bibitem {Petrich} Petrich, M.: Lectures in Semigroups. Akademie-Verlag Berlin (1977)
\bibitem{Raney} Raney, G.N.: Completely Distributive Complete Lattices. Proceedings of the American Mathematical Society
3-5, 677--680 (1952)
\bibitem{Rosenfeld} Rosenfeld, A.: Fuzzy groups. Journal of Mathematical Analysis and Applications 35, 512--517 (1971)
\bibitem{Schein} Schein, B.M: Homomorphisms and subdirect decompositions of semigroups. Pacific Journal of Mathematics 17-3, 529--547 (1966)
\bibitem{Szasz} Sz\'asz, G.: Introduction to Lattice Theory. Akad\'emiai Kiad\'o, Budapest (1963)
\bibitem{Yamada} Yamada, M.: A note on subdirect decompositions of idempotent semigroups. Proceedings of the Japan Academy 36, 411-414 (1960)
\bibitem{Zadeh} Zadeh, L.A.: Fuzzy sets. Information and Control 8, 338--353 (1965)
\end{thebibliography}
\end{document}